\documentclass[leqno]{amsart}
\usepackage[text={6in,9in}, centering]{geometry}
\usepackage{bbm}
\usepackage{amsmath,amsthm,amssymb,amsfonts,amsthm,amsxtra}
\usepackage{color}
\usepackage{enumerate}

\def\dblsum{\mathop{\sum \sum}\limits}
\def\dblint{\mathop{\int_{1/2}^{5/2} \int_{1/2}^{5/2}}\limits}

\newtheorem*{acknow}{Acknowledgements}
\newtheorem*{remark}{Remark}
\newtheorem*{remark1}{Remark}
\newtheorem*{secondmoment}{Theorem 1 (Second Moment)}
\newtheorem*{subconvexity}{Corollary 1 (Subconvexity)}
\newtheorem*{thrmSS}{Theorem 2 (Shifted Convolution Sums)}

\newtheorem*{Integralbounds}{Lemma 1}
\newtheorem*{Derivatives}{Lemma 2}
\newtheorem*{Jeval}{Lemma 3}
\newtheorem*{Petersson}{Lemma 4 \textnormal{(Petersson trace formula)}}
\newtheorem*{KloostermanLS}{Lemma 5 \textnormal{(\cite{KMV} Prop. 5.1, \cite{DFI2}, \cite{DFI2E}, \cite{DI})}}

\newtheorem*{Voronoi}{Lemma 6 \textnormal{(Voronoi summation, \cite{KMV} Theorem A.4)}}

\newtheorem*{delta-symbol}{Lemma 7}

\newtheorem*{truncation}{Lemma 8}

\newtheorem*{shiftedsums}{Lemma 9}

\newtheorem*{Lemma2012}{Lemma 10}

\begin{document}
\title{Level aspect subconvexity for Rankin-Selberg $L$-functions} 
\author[]{Roman Holowinsky}
\address[]{Department of Mathematics, The Ohio State University, 100 Math Tower, 231 West 18th Avenue, Columbus, OH 43210-1174, USA.}
\email{holowinsky.1@osu.edu}
\author[]{Ritabrata Munshi}   
\address[]{School of Mathematics, Tata Institute of Fundamental Research, 1 Dr. Homi Bhabha Road, Mumbai 400005, India.}
\email{rmunshi@math.tifr.res.in}

\today

\begin{abstract}
Let $M$ be a square-free integer and let $P$ be a prime not dividing $M$ such that $P \sim M^\eta$ with $0<\eta<2/21$. We prove subconvexity bounds for $L(\tfrac{1}{2}, f \otimes g)$ when $f$ and $g$ are two primitive holomorphic cusp forms of levels $P$ and $M$. These bounds are achieved through an unamplified second moment method.
\end{abstract}

\maketitle


\section{Introduction and statement of results} 
Several authors have recently been successful in implementing the amplification method in order to establish level aspect subconvexity results for Rankin-Selberg convolutions of two $\rm{GL}(2)$ forms when one form is fixed and the other form is varying.   For example, if $f$ is a Hecke cusp form of fixed level and $g$ is a Hecke cusp form of varying level $M$, then various bounds of the form
$$
L(\tfrac{1}{2}, f \otimes g)\ll_f M^{1/2-\delta}
$$
for some absolute positive constant $\delta$ have been shown by Kowalski-Michel-VanderKam \cite{KMV}, Michel \cite{M} and Harcos-Michel \cite{HaM}.  Furthermore, results for the Rankin-Selberg convolution of two independently varying forms have been established in the works of Michel-Ramakrishnan \cite{MR}, Feigon-Whitehouse \cite{FW} and Nelson \cite{N} in situations where positivity of the central $L$-values is known.  Of particular interest, yet seemingly out of reach by means of current technology, are level and spectral aspect subconvexity results for the Rankin-Selberg convolution of two $\rm{GL}(2)$ forms of same level (e.g. when the two forms are same). These $L$-values appear naturally in many areas of number theory and in particular, have important connections with quantum chaos and equidistribution problems. \\

Subconvexity bounds for an individual $L$-function are often the result of sufficient bounds for a weighted average over an appropriate family of $L$-functions.  In this note, we consider the subconvexity problem for the Rankin-Selberg convolution of two varying $\rm{GL}(2)$ forms with co-prime levels through the use of a second moment method.  With the $L$-function here being constructed from data associated with two independently varying forms,  one has a large collection of natural families to choose from. \\ 

The ideas presented here may be applied to other Rankin-Selberg convolutions constructed out of multiple independently varying forms.  This is the first installment of recent work by the two authors related to the subconvexity problem and its purpose is to demonstrate the existence of situations in which subconvexity may be established through a second moment average without amplification.\\
  
\begin{acknow}
\textnormal{We thank IAS Princeton for the wonderful working conditions in which many of the ideas for this collaboration were initially conceived.  We also thank MSRI Berkeley, MF Oberwolfach and TIFR Mumbai for providing the opportunity for further discussions. The first author is supported by the Sloan fellowship BR2011-083 and the NSF grant DMS-1068043.}
\end{acknow}


\subsection{Holomorphic cusp forms}
 Let $N>0$ be an integer and $k>0$ be an even integer.  We denote by 
$\mathcal{S}_k(N)$ 
the linear space of holomorphic cusp forms of weight $k$, level $N$ and trivial nebentypus. Such forms are holomorphic functions on the upper half-plane $f:\mathbbm{H} \rightarrow \mathbbm{C}$ satisfying
\begin{equation}
f(\gamma z) = (cz+d)^k f(z) \label{modularity}
\end{equation}
for every $\gamma=\left(\begin{smallmatrix}a & b \\ c & d \end{smallmatrix}\right)\in \Gamma_0(N)$ and which vanish at every cusp. Any form $f\in \mathcal{S}_k(N)$ has a Fourier series expansion
$$
f(z)=\sum_{n\geqslant 1} \psi_f(n)n^{\frac{k-1}{2}} e(nz)
$$
with coefficients $\psi_f(n)$ satisfying 
$$
\psi_f(n) \ll_f \tau(n)
$$
as proven by Deligne \cite{D}.\\

The space $\mathcal{S}_k(N)$ is a finite dimensional Hilbert space with respect to the Petersson inner product
\begin{equation}
 \langle f_1, f_2 \rangle = \int_{\Gamma_0(N)\backslash\mathbbm{H}}y^{k} f_1(z) \bar{f}_2(z) \frac{dx dy}{y^2}.
\end{equation}
We can choose an orthogonal basis $\mathcal{B}_k(N)$ for $\mathcal{S}_k(N)$ which consists of common eigenfunctions of all the Hecke operators $T_n$ with $(n,N)=1$. That is, each $f\in \mathcal{B}_k(N)$ satisfies
$$
(T_n f )(z) = \frac{1}{\sqrt{n}} \sum_{\substack{ad=n\\ (a,N)=1}}  \left(\frac{a}{d}\right)^{k/2} \sum_{b \, (\textnormal{mod }d)} f\left(\frac{az+b}{d}\right)= \lambda_f(n) f(z)
$$
for all $(n,N)=1$. Such $f$ are called \textit{Hecke eigen cusp forms}.  The Hecke operators are multiplicative and one has that 
$$
\psi_f(m)\lambda_f(n)=\sum_{d|(m,n)} \psi_f\left(\frac{mn}{d^2}\right)
$$
for any $m,n \geqslant 1$ with $(n,N)=1$. In particular, $\psi_f(1)\lambda_f(n) = \psi_f(n)$ if $(n,N)=1$.  Therefore,
\begin{equation}
\lambda_f(m) \lambda_f(n) = \sum_{d|(m,n)} \lambda_f\left(\frac{mn}{d^2}\right) \label{Heckerelations}
\end{equation}
if $(nm,N)=1$.   The Hecke eigenbasis $\mathcal{B}_k(N)$ also contains a subset of \textit{newforms} $\mathcal{B}^\ast_k(N)$, those forms which are simultaneous eigenfunctions of all the Hecke operators $T_n$ for any $n\geqslant 1$ and normalized to have first Fourier coefficient $\psi_f(1)=1$. For $f\in \mathcal{B}^\ast_k(N)$, the Hecke relations \eqref{Heckerelations} hold for all integers $n,m \geqslant 1$ and it is also known (see \cite{ILS}) that 
\begin{eqnarray}
|\lambda_f(p)| = p^{-1/2} & \textnormal{for any} & p|N. \label{Fricke}
\end{eqnarray}\\


\subsection{Rankin-Selberg convolutions of forms with co-prime levels}
Let $N$ and $M$ be two positive square-free co-prime integers and let $k$ and $\kappa$ be two fixed positive even integers.  Given two newforms $f\in \mathcal{B}^{\ast}_k(N)$ and $g\in \mathcal{B}^{\ast}_\kappa(M)$, we consider the associated Rankin-Selberg convolution $L$-function (see \cite{HaM})
$$
L(s, f \otimes g) = \prod_p \prod_{i=1}^2\prod_{j=1}^2\left(1-\frac{\alpha_{f,i}(p)\alpha_{g,j}(p)}{p^s}\right)^{-1} = \zeta^{(NM)}(2s)\sum_{n\geqslant 1} \lambda_f(n) \lambda_g(n) n^{-s}
$$
where the $\{\alpha_{f,i}\}$ and $\{\alpha_{g,j}\}$ are the local parameters of the $L$-functions associated to $f$ and $g$ respectively and $\zeta^{(NM)}(2s)$ is the partial Riemann zeta function with the local factors at primes dividing $NM$ removed. The local parameters satisfy the relations $\alpha_{f,1}(p)+\alpha_{f,2}(p)=\lambda_f(p)$ and $\alpha_{f,1}(p)\alpha_{f,2}(p)=\chi_0(p)$ with $\chi_0$ the principal character of modulus $N$ and similarly for the local parameters associated with $g$. The completed $L$-function is then defined as
$$
\Lambda(s, f\otimes g):=\mathcal{Q}^{s/2} L_\infty(s, f\otimes g,s) L(s, f\otimes g)
$$
where the conductor (see \cite{M}) is given by
$\mathcal{Q}:=\mathcal{Q}(f\otimes g)=(NM)^2$
and the local factor at infinity (see \cite{IK}) is a product of gamma factors
$$
L_\infty(s, f \otimes g):=\pi^{-2s} \Gamma\bigg(\frac{s+\frac{|k-\kappa|}{2}}{2}\bigg)\Gamma\bigg(\frac{s+\frac{k+\kappa}{2}}{2}\bigg)\Gamma\bigg(\frac{s+\frac{|k-\kappa|}{2}+1}{2}\bigg)\Gamma\bigg(\frac{s+\frac{k+\kappa}{2}-1}{2}\bigg).
$$
The completed $L$-function satisfies the functional equation
$$
\Lambda(s, f\otimes g)= \Lambda(1-s, f\otimes g).
$$\\

\begin{remark1}
\textnormal{We have restricted our discussion to the case of trivial nebentypus as historically this has been the hardest case (see \cite{HaM}, \cite{KMV}, \cite{M}).  We have taken $(N,M)=1$ to ensure that the conductor is as large as possible.  For general $N$ and $M$ we have that $(NM)^2/(N,M)^4\leqslant\mathcal{Q}(f\otimes g)\leqslant (NM)^2/(N,M)$ (see \cite{HaM}).}\\
\end{remark1}

The convexity bound for  $L(s, f\otimes g)$ at the point $s=1/2$ is
$$
L(\tfrac{1}{2}, f\otimes g) \ll_\varepsilon \mathcal{Q}^{1/4+\varepsilon}
$$
for any $\varepsilon>0$ and may be established in this case simply by the approximate functional equation and Deligne's bound. It has recently been shown by Heath-Brown \cite{H2}, in the general setting of Selberg class $L$-functions using Jensen's formula for strips, that the $\varepsilon$ in the above bound may be removed
 $$
L(\tfrac{1}{2}, f\otimes g) \ll \mathcal{Q}^{1/4}.
$$
Furthermore, the general results of Soundararajan \cite{S} provide a ``weak-subconvexity'' bound of the form
$$
L(\tfrac{1}{2}, f \otimes g) \ll \frac{\mathcal{Q}^{1/4}}{(\log \mathcal{Q})^{1-\varepsilon}}
$$  
for any $\varepsilon>0$.  \\

\subsection{Main results}
Our purpose here is to provide level aspect subconvexity bounds for the Rankin-Selberg convolution of two forms of varying levels $N$ and $M$ in situations where both forms are varying at different rates, say $N\sim M^{\eta}$ for some $0<\eta<1$.  The main point we wish to stress, is that we take advantage of the size of the smaller level $N$.  The method we present here does not produce subconvexity bounds when $N=1$ nor when $N$ is the same size as $M$.  Both levels must contribute to the complexity of the problem and they must do so in a manner which is sufficiently distinguishable for the method to work.  We restrict to the case of $N=P$ prime to simplify our presentation. Recall that our conductor in this case is of size $\mathcal{Q}=(PM)^2$.\\

We start by reducing our $L$-function to a smooth sum over Hecke eigenvalues by a standard approximate functional equation argument, see for example \cite{IK}, \cite{IM}, \cite{M}. Since we are working with newforms of trivial nebentypus, we have
$$
L(\tfrac{1}{2}, f\otimes g) = 2\sum^{\infty}_{n=1}\frac{\lambda_f(n)\lambda_g(n)}{\sqrt{n}} W\left(\frac{n}{\sqrt{\mathcal{Q}}}\right)
$$
where 
$$
W(y)=\frac{1}{2\pi i} \int_{(3)} G(u) \frac{L_\infty(\tfrac{1}{2}+u, f\otimes g)}{L_\infty(\tfrac{1}{2}, f\otimes g)} \zeta^{(NM)}(1+2u) y^{-u} \frac{du}{u}
$$
and
$$
G(u)=\left(\cos\frac{\pi u}{4A}\right)^{-16A}
$$
for any positive integer $A$.  The derivatives of $W(y)$ satisfy
$$
y^j W^{(j)}(y)\ll_{k,\kappa} \mathcal{Q}^{\varepsilon}(1+y)^{-A} \log(2+y^{-1})
$$
for any $\varepsilon>0$. Applying a smooth partition of unity one may derive that (see e.g. \cite{IM})
$$
L(\tfrac{1}{2}, f\otimes g)\ll_{k,\kappa} \mathcal{Q}^{\varepsilon} \sum_{X} \frac{|L_{f\otimes g}(X)|}{\sqrt{X}}\left(1+\frac{X}{\sqrt{\mathcal{Q}}}\right)^{-A}
$$
where
$$
L_{f\otimes g}(X)=\sum_n \lambda_f(n)\lambda_g(n) h\left(\frac{n}{X}\right)
$$
and $h$ is a smooth function, compactly supported on $[\frac{1}{2},\frac{5}{2}]$ with bounded derivatives and $X$ runs over values $2^\nu$ with $\nu=-1,0,1,2,\ldots$.  \\

Since  $L_{f\otimes g}(X)$ is trivially bounded by $X^{1+\varepsilon}$ for any $\varepsilon>0$, the contribution from those $X>\mathcal{Q}^{1/2+\varepsilon}$ is made negligible by choosing $A$ above to be sufficiently large.  Likewise, if $X<\mathcal{Q}^{1/2-\delta}$ for some $\delta>0$, then  $L_{f\otimes g}(X)X^{-1/2} \ll \mathcal{Q}^{1/4-\delta/2}$. Therefore, we are left with
$$
L(\tfrac{1}{2}, f\otimes g)\ll_{\varepsilon} \mathcal{Q}^{\varepsilon}\left\{ \mathcal{Q}^{1/4-\delta/2}+\max_{\mathcal{Q}^{1/2-\delta} \leqslant X\leqslant \mathcal{Q}^{1/2+\varepsilon}} 
\frac{|L_{f\otimes g}(X)|}{\sqrt{X}}\right\}
$$
for any $\delta>0$. Subconvexity bounds will now follow if one is able to sufficiently bound $L_{f\otimes g}(X)$ in the remaining range for $X$.  We shall do so by averaging over a Hecke eigenbasis for forms of level $M$.\\

\begin{secondmoment}
Let $M$ be a positive square-free integer and let $P$ be a prime such that $(P,M)=1$. Let $k$ and $\kappa$ be two fixed positive even integers.  Set $\mathcal{Q}=(PM)^2$.  Let $\varepsilon,\delta>0$ and choose any $\mathcal{Q}^{1/2-\delta} \leqslant X\leqslant \mathcal{Q}^{1/2+\varepsilon}$.  For any new form $f \in \mathcal{B}^{\ast}_k(P)$ we have
$$
\sum_{g\in \mathcal{B}_\kappa(M)} \omega_g^{-1} \left|\sum_n \psi_f(n)\psi_g(n) h\left(\frac{n}{X}\right)\right|^2 \ll_{\varepsilon,\delta} X P\mathcal{Q}^{\varepsilon} \left(\frac{1}{P}+\frac{1}{\mathcal Q^\delta}+\mathcal{Q}^{\frac{5}{4}\delta}\frac{P^{\frac{21}{8}}}{M^{\frac{1}{4}}}+\mathcal{Q}^{3\delta}\frac{P^2}{M^{\frac{1}{4}}}\right)
$$ 
where the spectral weights are given as $\omega_g:=\frac{(4\pi)^{\kappa-1}}{\Gamma(\kappa-1)}\langle g, g \rangle$.\\
\end{secondmoment}

Note that a second moment bound of the form 
\begin{equation}
\sum_{g\in \mathcal{B}_\kappa(M)} \omega_g^{-1}\left|\sum_n \psi_f(n)\psi_g(n) h\left(\frac{n}{X}\right)\right|^2 \ll X P \mathcal{Q}^{\varepsilon}\label{secondmomentconvexity}
\end{equation}
for all $X\leqslant \mathcal{Q}^{1/2+\varepsilon}$ and any $\varepsilon>0$ would produce the convexity bound for any individual $L(\tfrac{1}{2}, f \otimes g)$ with $f$ and $g$ both newforms since then $\psi=\lambda$ and $\omega_g \ll_\kappa M$ (see \cite{ILS}).  Therefore, the bound in Theorem 1 produces a subconvexity bound when $P \sim M^\eta$ with $0<\eta<2/21$. \\

\begin{subconvexity}
Let $M$ be a positive square-free integer and let $P$ be a prime not dividing $M$. Let $\eta=\frac{\log P}{\log M}$. Let $k$ and $\kappa$ be two fixed positive even integers.  For two newforms $f \in \mathcal{B}^{\ast}_k(P)$ and $g \in \mathcal{B}^{\ast}_\kappa(M)$ we have
$$
L(\tfrac{1}{2}, f \otimes g)\ll \mathcal{Q}^{\frac{1}{4}+\varepsilon} \left(\frac{1}{\mathcal{Q}^{\frac{\eta}{2(1+\eta)}}}+\frac{1}{\mathcal Q^{\frac{2-21\eta}{64(1+\eta)}}}\right).
$$
\begin{proof}
Soften the bound in Theorem 1 to
$$
\sum_{g\in \mathcal{B}_\kappa(M)} \omega_g^{-1} \left|\sum_n \psi_f(n)\psi_g(n) h\left(\frac{n}{X}\right)\right|^2 \ll_{\varepsilon,\delta} X P\mathcal{Q}^{\varepsilon} \left(\frac{1}{P}+\frac{1}{\mathcal Q^\delta}+\mathcal{Q}^{3\delta}\frac{P^{\frac{21}{8}}}{M^{\frac{1}{4}}}\right)
$$
and equate the second and third terms on the right hand side above while replacing all occurrences of $P$ by $M^{\eta}$.
\end{proof}
\end{subconvexity}

The estimates that we have obtained in Theorem 1 and Corollary 1 are the result of analysis of the shifted convolution sum problem through the $\delta$-method (\cite{DFI}, \cite{H}) with explicit dependence on the level $P$ of the form $f$.  It is possible to push our arguments further to improve these estimates by considering the shifted convolution sum problem on average over shifts while again maintaining explicit dependence on the level $P$ of $f$ and we shall do so in a later work.  For our purposes here, we prove the following theorem for a fixed non-zero shift.\\

\begin{thrmSS}
Let $\ell$ be a non-zero integer and let $X,Y\geqslant 1$. Let $F$ be a smooth function supported on $[1/2,5/2]\times[1/2,5/2]$ with partial derivatives satisfying
$$
x^iy^j \frac{\partial^i}{\partial x^i}\frac{\partial^j}{\partial y^j} F\left(\frac{x}{X},\frac{y}{Y}\right) \ll Z {Z_x}^{i}{Z_y}^{j}
$$
for some $Z>0$ and $Z_x,Z_y\geqslant 1$.  For any new forms $f_1, f_2 \in\mathcal{B}^\ast_k(P)$ we have
$$
\dblsum_{m=nP+\ell} \lambda_{f_1}(n)\lambda_{f_2}(m)F\left(\frac{n}{X},\frac{m}{Y}\right) \ll
(XYP)^{\varepsilon}P \max\{XP,Y\}^{3/4}Z\sqrt{Z_xZ_y}\max\{Z_x,Z_y\}^{5/4}.
$$
\end{thrmSS}

For other works involving estimates of shifted sums see \cite{B}, \cite{BH}, \cite{DFI}, \cite{DFI-2}, \cite{HaM}, \cite{Ho}, \cite{Ho2}, \cite{J}, \cite{J2}, \cite{LS}, \cite{P}, \cite{R}, \cite{Sa} and \cite{Ha} for dependence on the level of the forms.  The above bound in Theorem~2 does not follow easily from any of the above works. The main advantage here is uniformity with respect to the shift $\ell$ and the coefficient $P$.  Furthermore, we note that if $\ell\equiv 0 \bmod{P}$ then one also has the trivial bound $ZX/\sqrt{P}$ by using \eqref{Fricke}.\\


\section{Preliminaries}

\subsection{Bessel functions}

We record here some standard facts about the $J$-Bessel functions as can be seen in \cite{W} as well as several estimates for integrals involving Bessel functions which will be required for our application.  One may write the $J$-Bessel functions as 
\begin{equation}
J_k(x) = e^{ix} W_k(x)+e^{-ix} \overline{W}_k(x) \label{Joss}
\end{equation}
where 
\begin{equation}
W_k(x)=\frac{e^{i(\frac{\pi}{2}k-\frac{\pi}{4})}}{\Gamma(k+\frac{1}{2})}\sqrt{\frac{2}{\pi x}} \int_0^{\infty} e^{-y}(y(1+\frac{iy}{2x}))^{k-\frac{1}{2}}dy \label{JW}
\end{equation}
which, when $k$ is a positive integer, one has that
\begin{equation}
x^j W_k^{(j)}(x) \ll \frac{x}{(1+x)^{3/2}}. \label{JWbounds}
\end{equation}
Using the above facts leads us to the following results.\\

\begin{Integralbounds}
Let $k, \kappa\geqslant 2$ be integers and let $a,b,x,y>0$. Define
$$
I(x,y):=\int_0^{\infty} h\left(\xi\right) J_{\kappa-1}\left(4\pi a\sqrt{x\xi}\right)J_{k-1}\left(4\pi b\sqrt{y\xi}\right)d\xi
$$
where $h$ is a smooth function compactly supported on $\left[\frac{1}{2},\frac{5}{2}\right]$ with bounded derivatives. We have 
$$
I(x,y) \ll_j |a\sqrt{x}-b\sqrt{y}|^{-j}
$$
for any $j\geqslant 0$.
\begin{proof}
A change of variables, $\xi=w^2$, gives
$$
I(x,y)=2\int_0^{\infty} h(w^2)\; w\;  J_{\kappa-1}\left(4\pi a\sqrt{x} w\right)J_{k-1}\left(4 \pi b\sqrt{y} w\right)dw.
$$
Therefore, we see from \eqref{Joss} that $I(x,y)$ may be written as the sum of four similar terms, one of them being
$$
\int_0^{\infty} e\left(2w(a\sqrt{x}-b\sqrt{y}))\right)h(w^2)\;w \;W_{\kappa-1}\left(4\pi a\sqrt{x}w\right)\overline{W}_{k-1}\left(4\pi b\sqrt{y}w\right)dw.
$$
Repeated integration by parts gives the desired result.
\end{proof}
\end{Integralbounds}

\vspace{.5cm}

\begin{Derivatives}
For $I(x,y)$ as in Lemma 1, we have
$$
x^iy^j \frac{\partial^i}{\partial x^i}\frac{\partial^j}{\partial y^j} I(x,y) \ll_{i,j}   \frac{a\sqrt{x}}{(1+a\sqrt{x})^{3/2}} \frac{b\sqrt{y}}{(1+b\sqrt{y})^{3/2}}\left(1+a\sqrt{x}\right)^{i}\left(1+b\sqrt{y}\right)^{j}. 
$$
\begin{proof}
Differentiate and use the bound in \eqref{JWbounds}.
\end{proof}
\end{Derivatives}

\begin{Jeval}
Let $k,P,q$ be positive integers with $k\geqslant 2$ and let $\ell$ be a non-zero integer. Take $Q>1$ and $X,Y\geqslant 1$. For any $a,b>0$, define
\begin{equation}
J(a,b) := \int_0^{\infty}\int_0^{\infty}F\left(\frac{x}{X},\frac{y}{Y}\right) h\left(\frac{q}{Q},\frac{xP+\ell-y} {Q^2}\right)J_{k-1}\left(4 \pi a\sqrt{x}\right)J_{k-1}\left(4\pi b \sqrt{y}\right)dx\; dy
\end{equation}
where $h\left(\frac{q}{Q},\frac{xP+\ell-y} {Q^2}\right)$ is the function from Lemma~7 in \S2.2 and $F$ is a smooth function supported on $[1/2,5/2]\times[1/2,5/2]$ with partial derivatives satisfying
$$
x^iy^j \frac{\partial^i}{\partial x^i}\frac{\partial^j}{\partial y^j} F\left(\frac{x}{X},\frac{y}{Y}\right) \ll Z {Z_x}^{i}{Z_y}^{j}
$$
for some $Z>0$ and $Z_x,Z_y\geqslant 1$.  We have
\begin{equation}
J(a,b) \ll  ZXY\frac{Q}{q} \frac{a\sqrt{X}}{(1+a\sqrt{X})^{3/2}} \frac{b\sqrt{Y}}{(1+b\sqrt{Y})^{3/2}}\left[\frac{1}{a\sqrt{X}}\left\{Z_x+\frac{XP}{qQ}\right\}\right]^{i}\left[\frac{1}{b\sqrt{Y}}\left\{Z_y+\frac{Y}{qQ}\right\}\right]^{j}\label{Jbounds}
\end{equation}
for any non-negative integers $i$ and $j$.  Furthermore,
\begin{equation}
J(a,b) \ll  \frac{ZXY}{(1+a\sqrt{X})^{3/2} \; (1+b\sqrt{Y})^{3/2}} \; \frac{Q}{q} \; \min\{Z_x \; b \sqrt{Y}, Z_y \; a \sqrt{X}\}\; Q^{\varepsilon}.\label{JsaveXbyqQ}
\end{equation}\\
\begin{proof}
A change of variables, integrating by parts once in $x$ and applying the given bounds for the functions $F, h$ and the Bessel functions gives
$$
J(a,b) \ll ZXY\frac{Q}{q} \frac{a\sqrt{X}}{(1+a\sqrt{X})^{3/2}} \frac{b\sqrt{Y}}{(1+b\sqrt{Y})^{3/2}}\left[\frac{1}{a\sqrt{X}}\left\{Z_x+XPI\right\}\right]
$$
with
$$
I:= \dblint_{2|xXP+\ell-yY|>qQ} \frac{1}{|xXP+\ell-yY|}\;dx\;dy.
$$
Trivially, $I \ll (qQ)^{-1}$ and this is how one arrives at \eqref{Jbounds} with $i=1$ and $j=0$. Repeated integration by parts would then establish \eqref{Jbounds} for all $i$ and $j$.   Otherwise, replace $x$ by $u=xXP+\ell - yY$ so that $dx=(XP)^{-1} du$ and
$$
I\ll (XP)^{-1} \int_{1/2}^{5/2} \int_{qQ/2}^{(XP+Y+|\ell|)Q^{\varepsilon}} \frac{1}{u} \;du \; dy \ll (XP)^{-1} Q^{\varepsilon}.
$$
Repeating the argument, for $y$ instead of $x$, gives the bound \eqref{JsaveXbyqQ}.
\end{proof}
\end{Jeval}

\subsection{Summation Formulae, Large Sieve and the $\delta$-method}
Let $k\geqslant 2$ be an integer.  For any $n,m,c \in \mathbbm{N}$, let $S(n,m;c)$ denote the Kloosterman sum
$$
S(n,m;c)=\sideset{}{^\ast}\sum_{\alpha(c)} e\left(\frac{n \alpha+m\overline{\alpha}}{c}\right).
$$ 
The Kloosterman sums satisfy the Weil bound
$$
|S(n,m;c)| \leqslant (n,m,c)^{1/2} c^{1/2} \tau(c)
$$
where $\tau(c)$ is the number of divisors of $c$. This bound is best possible for an individual Kloosterman sum.  Sums of Kloosterman sums appear in the following spectral average (see \cite{IK} for a derivation).\\
\begin{Petersson} Let $N\geqslant 1$ be an integer. Let $\mathcal{B}_k(N)$ be any Hecke eigenbasis for $\mathcal{S}_k(N)$. For any $n,m\geqslant 1$, we have
$$
 \sum_{f\in \mathcal{B}_k(N)} \omega_f^{-1}\psi_f(n)\overline{\psi_f(m)} = \delta(n,m)+ 2\pi i^{-k} \sum_{\substack{c>0\\c\equiv 0 (N)}}\frac{1}{c} S(n,m;c) J_{k-1}\left(\frac{4\pi \sqrt{nm}}{c}\right)
$$
where the spectral weights $\omega_f$ are given by
$$
\omega_f:= \frac{(4\pi)^{k-1}}{\Gamma(k-1)}\langle f, f \rangle
$$
and $\delta(n,m)=1$ if $n=m$ and $\delta(n,m)=0$ otherwise.
\end{Petersson}

\noindent One also has the following large sieve estimate.\\

\begin{KloostermanLS}
Let $\eta$ be a smooth function supported on $[C/2,5C/2]$ such that $\eta^{(j)} \ll_j C^{-j}$ for all $j\geqslant 0$.  For any sequences of complex numbers $x_n,y_m$ we have
\begin{align*}
\sum_{n\leqslant X} \sum_{m\leqslant Y} x_n y_m \sum_{\substack{c>0\\ c\equiv 0 (N)}}&\frac{\eta(c)}{c} S(n,m;c) J_{k-1}\left(\frac{4\pi \sqrt{nm}}{c}\right)\\
&\ll_{\varepsilon,k} C^{\varepsilon} \bigg(\frac{\sqrt{XY}}{C}\bigg)^{k-3/2} \bigg(1+\frac{X}{N}\bigg)^{1/2}\bigg(1+\frac{Y}{N}\bigg)^{1/2} \|x\|_2 \|y\|_2
\end{align*}
with any $\varepsilon>0$. Moreover the exponent $k-3/2$ may be replaced by $1/2$.\\
\end{KloostermanLS}

The above estimate will be useful in controlling the size of Kloosterman sum moduli. For all remaining moduli we will apply the following analogue to Poisson summation.\\

\begin{Voronoi}
Let $(a,q)=1$ and let $h$ be a smooth function, compactly supported in $(0,\infty)$. Let $f$ be a holomorphic newform of level $N$ and weight $k$. Set $N_2:=N/(N,q)$. Then there exists a complex number $\eta$ of modulus $1$ (depending on $a,q$ and $f$) and a newform $f^\ast$ of the same level $N$ and the same weight $k$ such that
$$
\sum_n \lambda_f(n) e\bigg(n\frac{a}{q}\bigg)h(n) =  \frac{2 \pi\eta}{q\sqrt{N_2}}\sum_n \lambda_{f^{\ast}}(n) e\bigg(-n\frac{\overline{aN_2}}{q}\bigg)\int_0^{\infty} h(\xi) J_{k-1}\bigg(\frac{4\pi\sqrt{n\xi}}{q\sqrt{N_2}}\bigg) d\xi
$$
where $\overline{x}$ denotes the multiplicative inverse of $x$.
\end{Voronoi}

We will now briefly recall a version of the circle method introduced in \cite{DFI} and \cite{H}. The starting point is a smooth approximation of the $\delta$-symbol. We will follow the exposition of Heath-Brown in \cite{H}.\\

\begin{delta-symbol}
For any $Q>1$ there is a positive constant $c_Q$, and a smooth function $h(x,y)$ defined on $(0,\infty)\times\mathbb R$, such that
\begin{align}
\label{cm}
\delta(n,0)=\frac{c_Q}{Q^2}\sum_{q=1}^{\infty}\;\sideset{}{^{\star}}\sum_{a \bmod{q}}e\left(\frac{an}{q}\right)h\left(\frac{q}{Q},\frac{n} {Q^2}\right).
\end{align}
The constant $c_Q$ satisfies $c_Q=1+O_A(Q^{-A})$ for any $A>0$. Moreover $h(x,y)\ll x^{-1}$ for all $y$, and $h(x,y)$ is non-zero only for $x\leq\max\{1,2|y|\}$. 
\end{delta-symbol} 

In practice, to detect the equation $n=0$ for a sequence of integers in the range $[-X,X]$, it is logical to choose $Q=X^{1/2}$.  The smooth function $h(x,y)$ satisfies (see \cite{H}) 
\begin{eqnarray}
x^{i} \frac{\partial^i}{\partial x^i}h(x,y)\ll_i x^{-1} & \textnormal{and} & \frac{\partial}{\partial y}h(x,y)=0 \label{hbound1}
\end{eqnarray}
for $x\leq 1$ and $|y|\leq x/2$. Also for $|y|>x/2$, we have
\begin{equation}
x^i y^j \frac{\partial^{i}}{\partial x^i}\frac{\partial^{j}}{\partial y^j}h(x,y)\ll_{i,j} x^{-1}. \label{hbound2}
\end{equation}


\section{Initial reduction of the second moment}

Let $M$ be a positive square-free integer and let $P$ be a prime not dividing $M$.  Let $k$ and $\kappa$ be two positive even fixed integers. Fix a newform $f\in \mathcal{B}^{\ast}_k(P)$ and choose an orthogonal Hecke eigenbasis $\mathcal{B}_\kappa(M)$ for $\mathcal{S}_\kappa(M)$. Set $\mathcal{Q}:=(PM)^2$.  Let $\varepsilon,\delta>0$ and choose any $\mathcal{Q}^{1/2-\delta}\leqslant X\leqslant \mathcal{Q}^{1/2+\varepsilon}$.  As seen in the statement of Theorem 1, we are interested in obtaining upper bounds for the sum
\begin{equation}
S_f(X):=\sum_{g\in \mathcal{B}_\kappa(M)} \omega_g^{-1} \bigg|\sum_n \psi_f(n)\psi_g(n) h\left(\frac{n}{X}\right)\bigg|^2 \label{secondmoment}
\end{equation}
where $\omega_g=\frac{(4\pi)^{\kappa-1}}{\Gamma(\kappa-1)}\langle g, g \rangle$ and $h$ is smooth, compactly supported on $[\frac{1}{2},\frac{5}{2}]$ with bounded derivatives.  We start by opening the square and applying the Petersson trace formula in $g$. Since $f$ is a newform, we have $\psi_f(n)=\lambda_f(n)$ and so
\begin{align*}
S_f(X) =& \sum_n \lambda_f(n)^2 h\left(\frac{n}{X}\right)^2\\
&+2\pi i^{-\kappa}\sum_n 
\sum_m \lambda_f(n) h\left(\frac{n}{X}\right) \lambda_f(m)h\left(\frac{m}{X}\right)\sum_{\substack{d>0\\d \equiv 0 (M)}} \frac{S(n,m;d)}{d}J_{\kappa-1}\bigg(\frac{4 \pi \sqrt{nm}}{d}\bigg).
\end{align*}
The ``diagonal term'' satisfies
$$
\sum_n \lambda_f(n)^2 h\left(\frac{n}{X}\right)^2 \ll X\mathcal{Q}^{\varepsilon}
$$
for any $\varepsilon>0$. This is the first term seen in the bound in Theorem 1. We are now left with the ``off-diagonal'' terms
$$
\sum_n \sum_m \lambda_f(n) h\left(\frac{n}{X}\right) \lambda_f(m)h\left(\frac{m}{X}\right) \sum_{\substack{d>0\\d \equiv 0 (M)}} \frac{S(n,m;d)}{d}J_{\kappa-1}\bigg(\frac{4 \pi \sqrt{nm}}{d}\bigg).
$$\\

We start by truncating the sum over $d$.  By the Weil bound for individual Kloosterman sums and bounds for the Bessel functions in \S 2.1, there exist positive values $A$ and $B$ such that the sum over $d$ may be truncated to those $d \leqslant X^A$ up to an error term of size at most $X^{-B}M^{-1}$.  For the remaining sum over $d\leqslant X^A$, we introduce another smooth partition of unity and break the sum into dyadic segments of size $D$, as we did with our $n$-sum above, so that we are left with sums of type
\begin{equation}
R_{f,D}(X):=\sum_n \sum_m \lambda_f(n) h\left(\frac{n}{X}\right) \lambda_f(m)h\left(\frac{m}{X}\right) \sum_{\substack{d>0\\d \equiv 0 (M)}} \frac{S(n,m;d)}{d}J_{\kappa-1}\bigg(\frac{4 \pi \sqrt{nm}}{d}\bigg)\eta_D(d) \label{offdiagdyadic}
\end{equation}
where $\eta_D$ is a smooth function supported on $[D/2,5D/2]$.  Note that $D$ must be of size at least $M$ by the congruence condition.  Furthermore, an application of Lemma~5 shows that 
$$
R_{f,D}(X)\ll \left(\frac{X}{D}\right)^{k-3/2}\left(1+\frac{X}{M}\right)X\mathcal{Q}^{\varepsilon} 
$$
which is smaller than the bound in Theorem 1 as soon as $D>X\mathcal Q^{2\delta}$.  Therefore, bounding the second moment in \eqref{secondmoment} has reduced to the following statement.\\

\begin{truncation}
Let $\delta>0$.  For any $\mathcal{Q}^{1/2-\delta}\leqslant X\leqslant \mathcal{Q}^{1/2+\varepsilon}$ we have
\begin{equation}
S_f(X) \ll_{\varepsilon,\delta}  \mathcal{Q}^{\varepsilon}\left(X + PX\mathcal Q^{-\delta}+\sum_{M\leqslant D\leqslant X\mathcal Q^{2\delta}}R_{f,D}(X)\right) \label{initialbound}
\end{equation}
where $R_{f,D}(X)$ is given by \eqref{offdiagdyadic} above and  $D$ runs over dyadic values.
\end{truncation}

\begin{remark}
\textnormal{With additional work, one might also eliminate all $D<X\mathcal Q^{-\theta}$, for some $\theta>0$ depending on $\delta$, in order to improve the final range of sizes $P$ relative to $M$ for which subconvexity is achieved. To keep our presentation short, we shall only show how one may remove $D<\sqrt{P} M \mathcal{Q}^{-\delta}$ (see Lemma~10).}\\
\end{remark}

We emphasize here the significance of the level $P$ in our problem. Note that the first term $X \mathcal{Q}^{\varepsilon}$ in the above bound \eqref{initialbound}, which came from the diagonal term after applying the Petersson trace formula in $g$, beats the convexity bound for $S_f(X)$ by $P$. If $P$ were fixed, then Lemma~8 would already be insufficient for subconvexity. \\

\section{Reduction to Shifted Convolution Sums}  

Let $\delta>0$.  We now proceed with the analysis of $R_{f,D}(X)$, as defined by \eqref{offdiagdyadic} above, when $M \leqslant D\leqslant X\mathcal Q^{2\delta}$ with $\mathcal{Q}^{1/2-\delta}\leqslant  X\leqslant \mathcal{Q}^{1/2+\varepsilon}$.  Opening the Kloosterman sums and changing the order of summation, one is left to study
\begin{equation}
\sum_{\substack{d>0\\d \equiv 0 (M)}} \frac{\eta_D(d)}{d} \sideset{}{^\ast}\sum_{\beta(d)}\sum_n\lambda_f(n)e\left(n\frac{\beta}{d}\right)h\left(\frac{n}{X}\right) \sum_m \lambda_f(m)e\left(m\frac{\overline{\beta}}{d}\right)h\left(\frac{m}{X}\right)J_{\kappa-1}\bigg(\frac{4 \pi \sqrt{nm}}{d}\bigg).\label{beforeSCS}
\end{equation}
As in the works \cite{HaM}, \cite{KMV} and \cite{M}, an application of Voronoi summation in $m$ and the evaluation of the resulting Ramanujan sums will lead to a collection of shifted convolution sums.  Switching from Kloosterman sums to Ramanujan sums in such a manner was already seen in the work of Goldfeld \cite{G}.  Since the application of Voronoi summation will be for a newform $f$ of level $P$ and therefore depends on the divisibility of $d$ by $P$, we first break apart our $d$ sum as
$$
\sum_{LR=P}\sum_{\substack{d>0\\ (d,L)=1\\ d \equiv 0 (RM)}}\frac{\eta_D(d)}{d}\sideset{}{^\ast}\sum_{\beta(d)}\sum_n\lambda_f(n)e\left(n\frac{\beta}{d}\right)h\left(\frac{n}{X}\right) \sum_m \lambda_f(m)e\left(m\frac{\overline{\beta}}{d}\right)h\left(\frac{m}{X}\right)J_{\kappa-1}\bigg(\frac{4 \pi \sqrt{nm}}{d}\bigg).
$$\\

Voronoi summation in $m$ then gives that the inner sum, up to a constant, is equal to
$$
\frac{1}{d\sqrt{L}}\sum_m \lambda_{f^{\ast}}(m)e\left(-m\frac{\beta\overline{L}}{d}\right)\int_0^{\infty}h\left(\frac{\xi}{X}\right)J_{\kappa-1}\bigg(\frac{4 \pi \sqrt{n\xi}}{d}\bigg)J_{k-1}\left(\frac{4 \pi \sqrt{m\xi}}{d\sqrt{L}}\right)d\xi.
$$
This produces a Ramanujan sum over $\beta$ for each modulus $d$, which we write as
$$
\sideset{}{^\ast}\sum_{\beta(d)} e\left(\frac{\beta(nL-m)}{d}\right)=\sum_{bc=d}\mu(b) \sum_{\beta(c)}e\left(\frac{\beta(nL-m)}{c}\right).
$$
Summing over $\beta$ will now produce a congruence condition between $n$ and $m$ modulo $c$.  Thus, we have reduced \eqref{beforeSCS} to the following.

\begin{shiftedsums}
Let $\delta>0$ and let $R_{f,D}(X)$ be as in \eqref{offdiagdyadic} with $\mathcal{Q}^{1/2-\delta}\leqslant X\leqslant \mathcal{Q}^{1/2+\varepsilon}$.  For any $M\leqslant D\leqslant X\mathcal Q^{2\delta}$ we have
$$
R_{f,D}(X) \ll \sum_{LR=P} \frac{1}{\sqrt{L}} \sum_{\substack{d>0\\ (d,L)=1\\d\equiv 0 (RM)}}\frac{\eta_D(d)}{d}\sum_{bc=d}\frac{1}{b}\big|\Sigma_d(L;c)\big|
$$
with shifted convolution sums
$$
\Sigma_d(L;c)=\sum_n\sum_{m\equiv nL(c)}\lambda_f(n)\lambda_{f^{\ast}}(m)I_d(n,m)
$$
where
$$
I_d(n,m)= h\left(\frac{n}{X}\right)\int_0^{\infty}h\left(\frac{\xi}{X}\right)J_{\kappa-1}\bigg(\frac{4 \pi \sqrt{n\xi}}{d}\bigg)J_{k-1}\left( \frac{4\pi\sqrt{m\xi}}{d\sqrt{L}}\right)d\xi.
$$
\end{shiftedsums}
\noindent In the above, $I_d(n,m)$ determines the main contribution in the sum over $n$ and $m$ which occurs when $n\sim X$ and $m=nL+O(dL(1+d/X)\mathcal{Q}^{\varepsilon})$. The other ranges of summation are negligible as can be seen by Lemma~1.


\section{Proof of Theorem 1}
Theorem~1 will follow after an appropriate treatment of the shifted convolution sums $\Sigma_d(L;c)$ in Lemma~9. We break this apart into cases according to the value of $L$.


\subsection{Treatment of the shifted sums $\Sigma_d(1;c)$}
Since we are dealing with forms of level $P$ prime, we only have two types of shifted convolution sums to consider, those with $L=P$ and those with $L=1$.  In the latter case, the moduli $d$ must be of size at least $PM$ by the congruence condition. Applying Lemma~1 and the bound $I_d(n,m)\ll X \min\{1,X/d\}$ obtained from Lemma~2 one has that
$$
\Sigma_d(1;c) \ll \mathcal{Q}^{\varepsilon}\frac{X^2}{d} \dblsum_{\substack{n\sim X\\m=n+O(\frac{d^2}{X}\mathcal{Q}^{\varepsilon})\\ m\equiv n (c)}} 1 \ll \frac{X^3}{d}\left(1+\frac{d^2}{X c}\right)\mathcal{Q}^{\varepsilon}, 
$$
so that this contribution to bounding $R_{f,D}(X)$ is
\begin{equation}
\sum_{\substack{d>0\\ d\equiv 0 (PM)}}\frac{\eta_D(d)}{d}\sum_{bc=d}\frac{1}{b}\big|\Sigma_d(1;c)\big| \ll_\varepsilon  \frac{X^2}{PM} \mathcal{Q}^{\varepsilon} \ll_\varepsilon X \mathcal{Q}^{\varepsilon}, \label{5.1bound}
\end{equation}
which matches the first term in \eqref{initialbound}.  


\subsection{Treatment of the zero shift in $\Sigma_d(P; c)$}
We now examine the case of $L=P$ and the contribution of the sums
$$
\frac{1}{\sqrt{P}} \sum_{\substack{d>0\\(d,P)=1\\d \equiv 0(M)}}\frac{\eta_D(d)}{d} \sum_{bc=d} \frac{1}{b} \big|\Sigma_{d}(P;c)\big|
$$
to  $R_{f,D}(X)$.  We first treat the ``zero shift'' in the shifted sums $\Sigma_d(P; c)$, i.e. when $m=nP$.  One has
\begin{equation}
\frac{1}{\sqrt{P}} \sum_{\substack{d>0\\(d,P)=1\\d \equiv 0(M)}}\frac{\eta_D(d)}{d} \sum_{bc=d} \frac{1}{b} \big|\sum_n \lambda_f(n)\lambda_{f^\ast}(nP) I_d(n,nP)\big| \ll_\varepsilon  \frac{X^2}{PM} \mathcal{Q}^{\varepsilon} \ll_\varepsilon X \mathcal{Q}^{\varepsilon} \label{5.2bound}
\end{equation}
by using the fact that $|\lambda_{f^{\ast}}(nP)|=|\lambda_{f^{\ast}}(n)\lambda_{f^{\ast}}(P)|=|\lambda_{f^{\ast}}(n)| P^{-1/2}$ (using \eqref{Fricke}) and again the bound $I_d(n,nP) \ll X \min\{1,X/d\}$.  This also matches the first term in \eqref{initialbound}.  In fact, for the same reasons, one may also show that 
$$
\frac{1}{\sqrt{P}} \sum_{\substack{d>0\\(d,P)=1\\d \equiv 0(M)}}\frac{\eta_D(d)}{d} \sum_{bc=d} \frac{1}{b} \big|\sum_n\sum_{\substack{m\equiv nP(c)\\m \equiv 0 (P)}}\lambda_f(n)\lambda_{f^{\ast}}(m) I_d(n,m)\big|\ll_\varepsilon X \mathcal{Q}^{\varepsilon} .
$$
However, we will not use this fact in what follows.


\subsection{Treatment of the non-zero shifts in $\Sigma_d(P; c)$} Finally, we are left with the non-zero shifts 
$$
\sum_n\sum_{\substack{m\equiv nP(c)\\m\neq nP}}\lambda_f(n)\lambda_{f^{\ast}}(m)I_d(n,m).
$$
By Lemma 1, we need only consider those $m\equiv nP(c)$ with $n\sim X$ and $m=nP+O(dP(1+d/X)\mathcal{Q}^{\varepsilon})$.  Therefore, the congruence in the inner sums may be rewritten as an equation 
\begin{equation}
\sum_{0\neq |r| \ll \frac{dP}{c} (1+\frac{d}{X})\mathcal{Q}^{\varepsilon}} \; \dblsum_{\substack{m=nP+c r}} \lambda_f(n)\lambda_{f^{\ast}}(m)I_d(n,m) . \label{asshifts}
\end{equation}
We proceed by taking a smooth partition of unity for the sum over $m$ writing
$$
I_d(n,m)=: X \sum_{Y} F\left(\frac{n}{X},\frac{m}{Y}\right)
$$
where $Y$ runs over values $2^v$ with $v=-1,0,1,2,\ldots$ such that $m=nP+c r$ is soluble when $m\sim Y$ and $F$ is supported on $[1/2,5/2]\times[1/2,5/2]$.   Furthermore, by Lemma~2 and the support of $F$, one has that
\begin{equation}
x^iy^j \frac{\partial^i}{\partial x^i}\frac{\partial^j}{\partial y^j} F\left(\frac{x}{X},\frac{y}{Y}\right) \ll \frac{X/d}{(1+X/d)^{3/2}}\frac{\sqrt{XY/(d^2P)}}{(1+\sqrt{XY/(d^2P)})^{3/2}}\left(1+\frac{X}{d}\right)^{i}\left(1+\frac{\sqrt{XY}}{d\sqrt{P}}\right)^{j} \label{Fbounds}
\end{equation}
for any non-negative integers $i$ and $j$. Therefore, we may split apart the sums in \eqref{asshifts} as
$$
X  \sum_{Y}  \sum_{0\neq |r| \ll \frac{dP}{c} (1+\frac{d}{X})\mathcal{Q}^{\varepsilon}} \; \dblsum_{m=nP+ c r}\lambda_f(n) \lambda_{f^{\ast}}(m)F\left(\frac{n}{X},\frac{m}{Y}\right)
$$
which is bounded by 
\begin{equation}
X^2 \sum_{Y} \frac{X/d}{(1+X/d)^{3/2}}\frac{\sqrt{XY/(d^2P)}}{(1+\sqrt{XY/(d^2P)})^{3/2}}\frac{dP(1+d/X)}{c}\mathcal{Q}^{\varepsilon} \label{firstnonzero}
\end{equation}
through an application of \eqref{Fbounds} with $i=j=0$.
For general $X$ and $d$, this may be bounded by $$X^2 \frac{dP}{c} \mathcal{Q}^{\varepsilon}.$$  However, in the case of $d\ll X\mathcal Q^{-\delta}$, one has that $m=nP+c r$ is soluble only when $Y\sim XP$ so that \eqref{firstnonzero} then satisfies the stronger bound 
$$
X^2\frac{X/d}{(1+X/d)^2}\frac{dP}{c}\mathcal{Q}^{\varepsilon} \ll X \frac{d^2P}{c}\mathcal{Q}^{\varepsilon}.
$$
Therefore, one has the following Lemma.\\

\begin{Lemma2012}
Let $\delta>0$.  For any $M\leqslant D\leqslant X\mathcal Q^{2\delta}$ we have
\begin{equation}
R_{f,D}(X) \ll X P^{3/2} \mathcal{Q}^{\varepsilon}. \label{trivialRbound}
\end{equation}
Furthermore, if $D\ll X\mathcal Q^{-\delta}$ then
\begin{equation}
R_{f,D}(X) \ll X P \mathcal{Q}^{\varepsilon} \left(\frac{D}{\sqrt{P}M}\right). \label{specialRbound}
\end{equation}\\
\end{Lemma2012}

Since the bound for $R_{f,D}(X)$ in \eqref{specialRbound} is better than the convexity bound in \eqref{secondmomentconvexity} when $D< \sqrt{P}M\mathcal Q^{-\delta}$, we may restrict now  to the case of $\sqrt{P}M\mathcal Q^{-\delta}\leqslant D \leqslant X\mathcal Q^{2\delta}$. The remaining task is to show that one can improve on the bound \eqref{trivialRbound} by more than $\sqrt{P}$ when $D$ is of that size.  \\

For such values of $D$, an application of Theorem~2 to the shifted convolution sums
$$
S_{X,Y}(c r):=  \dblsum_{m=nP+c r} \lambda_f(n)\lambda_{f^{\ast}}(m)F\left(\frac{n}{X},\frac{m}{Y}\right)
$$
gives
$$
S_{X,Y}(c r)\ll \mathcal{Q}^{\varepsilon}P\max\{XP,Y\}^{3/4}Z \sqrt{Z_xZ_y}\max\{Z_x,Z_y\}^{5/4},
$$
where
$$
Z=\frac{X/d}{(1+X/d)^{3/2}}\frac{\sqrt{XY/(d^2P)}}{(1+\sqrt{XY/(d^2P)})^{3/2}},\;\;\;Z_x=\left(1+\frac{X}{d}\right)\;\;\;\text{and}\;\;\;Z_y=\left(1+\frac{\sqrt{XY}}{d\sqrt{P}}\right).
$$
Hence the contribution of these non-zero shifts to Lemma~9 is bounded by
\begin{align}
\mathcal{Q}^{\varepsilon} P^{3/2}\sum_{\substack{d>0\\(d,P)=1\\d \equiv 0(M)}}\eta_D(d)\left(1+\frac{X}{d}\right)\sum_Y \max\{XP,Y\}^{3/4} Z\sqrt{Z_xZ_y}\max\{Z_x,Z_y\}^{5/4}.\label{non-zerocontribution}
\end{align}
First consider $\sqrt{P}M\mathcal Q^{-\delta} \leqslant D<X$. In this case, we have that $Y\ll XP \mathcal{Q}^{\varepsilon}$ and \eqref{non-zerocontribution} reduces to
\begin{align}
\label{non-zerocontribution1}
\mathcal{Q}^{\varepsilon} XP^{3/2}(XP)^{3/4} \sum_{\substack{d>0\\(d,P)=1\\d \equiv 0(M)}}&\frac{\eta_D(d)}{d} \left(\frac{X}{d}\right)^{5/4}\ll \mathcal{Q}^{\varepsilon} (XP)\frac{X^2P^{\frac{5}{4}}}{MD^{\frac{5}{4}}}\\
\nonumber & \ll  \mathcal{Q}^{\frac{5}{4}\delta+\varepsilon} (XP)\frac{X^2P^{\frac{5}{8}}}{M^{\frac{9}{4}}}\ll \mathcal{Q}^{\frac{5}{4}\delta+\varepsilon} (XP)\frac{P^{\frac{21}{8}}}{M^{\frac{1}{4}}}. 
\end{align}
Next consider $X \leqslant D\leqslant X\mathcal Q^{2\delta}$. In this case, we have that $Y\ll D^2P\mathcal Q^{\varepsilon}/X$ and \eqref{non-zerocontribution} reduces to
\begin{align}
\label{non-zerocontribution2}
\mathcal{Q}^{\varepsilon} XP^{3/2}\left(\frac{D^2P}{X}\right)^{3/4} \sum_{\substack{d>0\\(d,P)=1\\d \equiv 0(M)}}&\frac{\eta_D(d)}{d}\ll  \mathcal{Q}^{\varepsilon} (XP)\frac{D^{\frac{3}{2}}P^{\frac{5}{4}}}{MX^{\frac{3}{4}}}\\
\nonumber &\ll \mathcal{Q}^{3\delta+\varepsilon} (XP)\frac{X^{\frac{3}{4}}P^{\frac{5}{4}}}{M}\ll \mathcal{Q}^{3\delta+\varepsilon} (XP)\frac{P^{2}}{M^{\frac{1}{4}}}.
\end{align}
Combining \eqref{non-zerocontribution1} and \eqref{non-zerocontribution2} with \eqref{5.1bound} and \eqref{5.2bound} in \S5.1, \S5.2 and inserting these bounds into Lemma~8, completes the proof of Theorem~1.


\section{Proof of Theorem 2}

Let $X,Y\geqslant 1$ and let $F$ be a smooth function supported on $[1/2,5/2]\times[1/2,5/2]$ with partial derivatives bounded by 
\begin{equation}
x^iy^j \frac{\partial^i}{\partial x^i}\frac{\partial^j}{\partial y^j} F\left(\frac{x}{X},\frac{y}{Y}\right) \ll Z {Z_x}^{i}{Z_y}^{j} \label{genFderiv}
\end{equation}
for some $Z>0$ and $Z_x,Z_y\geqslant 1$.  Let $P$ be a prime, and let $k$ be a fixed positive even integer.   For any $f_1, f_2 \in \mathcal{B}_k^\ast(P)$ we consider the shifted convolution sums
\begin{equation}
S_{X,Y}(\ell): = \dblsum_{m=nP+\ell} \lambda_{f_1}(n)\lambda_{f_2}(m)F\left(\frac{n}{X},\frac{m}{Y}\right) \label{SCS}
\end{equation}
with $\ell$ a fixed non-zero integer satisfying $|\ell| \leqslant 10 (XP+Y)$ such that the sum is non-trivial.  Detecting the equation $m=nP+\ell$ in \eqref{SCS} through an application of the $\delta$-method gives
\begin{align}
 \label{afterdelta}
S_{X,Y}(\ell)=\frac{1}{Q^2} \sum_{q=1}^{\infty}\;\sideset{}{^{\star}}\sum_{a(q)}& e\left(\frac{a\ell}{q}\right) \sum_n  \lambda_{f_1}(n) e\left(\frac{anP}{q}\right)\\
\nonumber &\times \sum_m \lambda_{f_2}(m) e\left(\frac{-am}{q}\right)F\left(\frac{n}{X},\frac{m}{Y}\right)h\left(\frac{q}{Q},\frac{nP+\ell-m} {Q^2}\right)
\end{align}
up to a negligible error term with the function $h$ as in Lemma~7.  As mentioned in \S 2.2 one expects to take $Q$ to be roughly of size $\max\{ \sqrt{XP}, \sqrt{Y} \}$.  \\

\begin{remark}
\textnormal{Consider the case of $X\sim PM$ and $Y\sim P^2M$.  Such is the situation in our subconvexity application if one initially takes $X\sim \mathcal{Q}^{1/2}$ and Kloosterman sum moduli of size $D\sim PM$ in order to focus on the transition range of the Bessel function. Taking moduli $q$ of size up to $Q=P\sqrt{M}$ may therefore be regarded as a reduction of size $M$ to the conductor of the $n$ and $m$ sums. One then returns to Kloosterman sums, of moduli $q$ rather than $d$, by further applications of Voronoi summation.} 
\end{remark}

\subsection{Voronoi summation in $m$}
We are now set to treat $S_{X,Y}(\ell)$ in the form seen in display \eqref{afterdelta}.  Since we will be applying Voronoi summation to our sums in $n$ and $m$, the resulting sums will depend on the divisibility of the moduli $q$ by powers of $P$. Indeed, an application of Voronoi summation to the $m$-sum gives, up to a constant factor,
$$
\frac{1}{q\sqrt{P_q}}\sum_m \lambda_{{f_2}^{\ast}}(m) e\left(\frac{\overline{aP_q}m}{q}\right)\int_0^{\infty}F\left(\frac{n}{X},\frac{y}{Y}\right) h\left(\frac{q}{Q},\frac{nP+\ell-y} {Q^2}\right)J_{k-1}\left(\frac{4 \pi \sqrt{my}}{q\sqrt{P_q}}\right)dy
$$
where $P_q=P/(P,q)$.  Therefore, \eqref{afterdelta} reduces to 
\begin{align}
\label{mVoronoi}
\frac{1}{Q^2}  \sum_{q=1}^{\infty} \frac{1}{q\sqrt{P_q}}\sum_n \sum_m & \lambda_{f_1}(n)\lambda_{{f_2}^{\ast}}(m)S(\ell+nP,m \overline{P_q};q)\\
\nonumber &\times \int_0^{\infty}F\left(\frac{n}{X},\frac{y}{Y}\right) h\left(\frac{q}{Q},\frac{nP+\ell-y} {Q^2}\right)J_{k-1}\left(\frac{4 \pi \sqrt{my}}{q\sqrt{P_q}}\right)dy. 
\end{align}
Although we have gained the Kloosterman sum structure, an application of the Weil bound here would still be insufficient for our goal.
\subsection{Voronoi summation in $n$}
Define
\begin{equation}
J_\alpha(n,m;q):=J\left(\frac{\sqrt{nP^\alpha}}{q\sqrt{P_q}},\frac{\sqrt{m}}{q\sqrt{P_q}}\right)\label{Jalpha}
\end{equation}
where $P^\alpha=(q,P^2)$, $P_q=P/(P,q)$ and $J\left(\frac{\sqrt{nP^\alpha}}{q\sqrt{P_q}},\frac{\sqrt{m}}{q\sqrt{P_q}}\right)$ is the function in Lemma~3. Opening the Kloosterman sum in \eqref{mVoronoi} and applying Voronoi summation to the $n$-sum gives, up to a constant factor,
\begin{equation}
\frac{1}{Q^2}  \sum_{q} \frac{\sqrt{P^\alpha}}{q^2P_q}  \sum_n\sum_m \lambda_{{f_1}^\ast}(n)\lambda_{{f_2}^\ast}(m) S_\alpha(n,m,\ell;q) J_\alpha(n,m;q)\label{mnVoronoi}
\end{equation}
where
\begin{equation*}
S_\alpha(n,m,\ell;q) = 
\begin{cases} 
 S(\ell,(mP-n)\overline{P^2};q)&\textnormal{if } \alpha=0,\\
 S(\ell \overline{P},(m-n)\overline{P};q/P) S(\ell \overline{q/P}, m \overline{q/P};P)&\textnormal{if } \alpha=1,\\
S(\ell,m-nP;q)&\textnormal{if } \alpha=2.
\end{cases}
\end{equation*}\\

\subsection{Application of Weil bound}
We now break apart the sums in \eqref{mnVoronoi} according to the size of $q$. First, we note that the bound \eqref{Jbounds} in Lemma~3 allows one to truncate the $n$ and $m$ sums to be of size
\begin{eqnarray}
n\leqslant T_1:=\frac{q^2 P_q}{P^{\alpha} X}\left(Z_x+\frac{XP}{qQ}\right)^2 (XYP)^{\varepsilon} & \textnormal{ and } & m \leqslant T_2:=\frac{q^2 P_q}{Y}\left(Z_y+\frac{Y}{qQ}\right)^2 (XYP)^{\varepsilon}. \label{nmtruncation}
\end{eqnarray}
When the parameters are such that either $T_1<1$ or $T_2<1$ in \eqref{nmtruncation}, then one has arbitrary saving in these situations. Otherwise,  we apply the bound \eqref{JsaveXbyqQ} from Lemma~3  to $J_\alpha(n,m;q)$ and the Weil bound for Kloosterman sums in order to bound \eqref{mnVoronoi} by
\begin{align*}
&(XYP)^{\varepsilon}\frac{ZXY}{Q}  \sum_{q\leqslant Q} \frac{\sqrt{P^\alpha}}{q^3P_q}\sum_{n\leqslant T_1} \sum_{m\leqslant T_2} (\ell,q)^{1/2} q^{1/2} \left(\frac{q^2 P_q}{\sqrt{nmP^\alpha XY}}\right)^{3/2}\frac{\min\{Z_y\sqrt{nP^\alpha X}, Z_x\sqrt{mY} \}}{q\sqrt{P_q}}\\
& \ll (XYP)^{\varepsilon}\frac{Z(XY)^{1/4}}{Q} \sum_{\delta|\ell} \delta^{1/2} \sum_{\substack{q\leqslant Q\\(q,\ell)=\delta}} \frac{(T_1T_2)^{1/4}}{q^{1/2}P^{\alpha/4}} \; \min\{Z_y \sqrt{T_1P^{\alpha}X}, Z_x \sqrt{T_2 Y}\}.
\end{align*}
Bounding the minimum by the geometric mean, and using \eqref{nmtruncation}, we get the bound
$$
(XYP)^{\varepsilon}\frac{ZP}{Q}\sqrt{Z_xZ_y}\sum_{\delta|\ell} \delta^{2} \sum_{\substack{q\leqslant \frac{Q}{\delta}\\(q,\ell)=1}} q^{3/2} \;\left(Z_x+\frac{XP}{q\delta Q}\right)\left(Z_y+\frac{Y}{q\delta Q}\right),
$$ 
which is dominated by
\begin{eqnarray}
(XYP)^{\varepsilon}ZP\sqrt{Z_xZ_y}Q^{3/2}\left(Z_x+\frac{XP}{Q^2}\right)\left(Z_y+\frac{Y}{Q^2}\right).\nonumber
\end{eqnarray}
We bound the last expression by
\begin{equation}
(XYP)^{\varepsilon}ZP\sqrt{Z_xZ_y}Q^{3/2}\left(\max\{Z_x,Z_y\}+\frac{\max\{XP,Y\}}{Q^2}\right)^2 \label{second}
\end{equation}
Choosing
$
Q=\left(\frac{\max\{XP,Y\}}{\max\{Z_x,Z_y\}}\right)^{1/2}
$
in \eqref{second} produces the final bound 
$$
S_{X,Y}(\ell) \ll (XYP)^{\varepsilon}ZP\sqrt{Z_xZ_y}\max\{XP,Y\}^{3/4}\max\{Z_x,Z_y\}^{5/4}.
$$
\bibliographystyle{plain}	
\bibliography{LevelAspect}		
\end{document}